\documentclass[12pt]{amsart}

\usepackage{enumitem}
\usepackage{psfrag}
\usepackage{graphicx}
\usepackage{pinlabel}
\usepackage{graphicx}
\usepackage{amsmath}
\usepackage{amssymb}
\usepackage{amscd}
\usepackage{autobreak}
\usepackage{picinpar}
\usepackage{times}
\usepackage{pb-diagram}
\usepackage{graphicx}
\usepackage{wrapfig}
\usepackage{xspace}
\usepackage{hyperref}
\usepackage{url}
\usepackage{caption}
\usepackage{subcaption}
\usepackage{fancyhdr}
\usepackage{overpic}
\usepackage{color}
\usepackage[square,comma,sort&compress,numbers]{natbib}
\usepackage{latexsym}
\usepackage{mathrsfs}
\usepackage{amsfonts}
\usepackage{hyperref}
\usepackage{amsthm}
\usepackage{appendix}
\usepackage{pdfsync}

\theoremstyle{definition}
\newtheorem{theorem}{Theorem}[section]
\newtheorem{lemma}[theorem]{Lemma}

\newtheorem{corollary}[theorem]{Corollary}
\newtheorem{definition}[theorem]{Definition}

\newtheorem{remark}[theorem]{Remark}
\newtheorem{conjecture}{Conjecture}
\newtheorem*{theorem*}{Theorem}
\def\qed{\hfill{Q.E.D.}\smallskip}

\begin{document}

\title{\bf Deformation of discrete conformal structures on surfaces}
\author{Xu Xu}

\date{\today}

\address{School of Mathematics and Statistics, Wuhan University, Wuhan, 430072, P.R.China}
 \email{xuxu2@whu.edu.cn}



\begin{abstract}
In \cite{G3}, Glickenstein introduced the discrete conformal structures on polyhedral surfaces
in an axiomatic approach from Riemannian geometry perspective.
It includes Thurston's circle packings, Bowers-Stephenson's inversive distance circle packings
and Luo's vertex scalings as special cases.
In this paper, we study the deformation of Glickenstein's discrete conformal structures by combinatorial curvature flows.
The combinatorial Ricci flow for Glickenstein's discrete conformal structures on triangulated surfaces \cite{ZGZLYG}
is a generalization of Chow-Luo's combinatorial Ricci flow for Thurston's circle packings \cite{CL}
and  Luo's combinatorial Yamabe flow for vertex scalings \cite{L1}.
We prove that the solution of the combinatorial Ricci flow for Glickenstein's discrete conformal structures on triangulated surfaces
can be uniquely extended.
Furthermore, under some necessary conditions,
we prove that the solution of the extended combinatorial Ricci flow on a triangulated surface
exists for all time and converges exponentially fast for any initial value.
We further introduce the combinatorial Calabi flow for Glickenstein's discrete conformal structures on triangulated surfaces
and study the basic properties of the flow.
These combinatorial curvature flows provide effective algorithms
for finding piecewise constant curvature metrics on surfaces with prescribed combinatorial curvatures.
\end{abstract}

\maketitle

%


\section{Introduction}

Discrete conformal structure on polyhedral manifolds is a discrete analogue of the conformal structure on Riemannian manifolds,
which assigns the discrete metrics defined on the edges by scalar functions defined on the vertices.
In the literature, there are lots of different types of discrete conformal structures that have been extensively studied on surfaces,
including Thurston's circle packings \cite{T1}, Bowers-Stephenson's inversive distance circle packings \cite{BS}, Luo's vertex scalings \cite{L1} and others.
However, most of these discrete conformal structures were invented and studied individually.
In \cite{G3}, Glickenstein pioneered an axiomatic approach to
the Euclidean discrete conformal structures on surfaces from Riemannian geometry perspective.
Following \cite{G3},
Glickenstein-Thomas \cite{GT} introduced the hyperbolic and spherical discrete conformal structures on polyhedral surfaces
in an axiomatic approach.
Furthermore, Glickenstein-Thomas \cite{GT} studied the classification of the discrete conformal structures
on polyhedral surfaces and sorted out important classes. Recently, Xu-Zheng \cite{XZ4} gave a full classification of Glickenstein's
discrete conformal structures on surfaces.
According to the classification, Glickenstein's discrete conformal structures
include different types of circle packings and Luo's vertex scalings on polyhedral surfaces as special cases
and generalizes them to a very general context.
On the other hand, motivated by Bobenko-Pinkall-Springborn's celebrated work \cite{BPS},
Zhang-Guo-Zeng-Luo-Yau-Gu \cite{ZGZLYG} constructed independently some subclasses of Glickenstein's discrete conformal structures on surfaces using 3-dimensional hyperbolic geometry.
Motivated by Thurston's work \cite{T1},
Guo-Luo \cite{GL} also constructed independently some subclasses of Glickenstein's discrete conformal structures using different cosine laws in hyperbolic geometry.
Comparing to the existing special types of discrete conformal structures on surfaces, the geometry of
Glickenstein's discrete conformal structures on surfaces are not extensively studied.
Here we mention some of the progresses on the study of Glickenstein's discrete conformal structures on surfaces, which are far from complete.
Glickenstein \cite{G3} and Glickenstein-Thomas \cite{GT} established the variational principles for Glickenstein's discrete conformal structures on surfaces.
Zhang-Guo-Zeng-Luo-Yau-Gu \cite{ZGZLYG} established independently the variational principles for some subclasses of Glickenstein's discrete conformal structures on surfaces and gives the explicit geometric interpolations of the action functional via 3-dimensional hyperbolic geometry.
Guo-Luo \cite{GL} also established independently the variational principles for some other subclasses of Glickenstein's discrete conformal structures on surfaces from de Verdi\`{e}re's viewpoint \cite{DV}, and studied the rigidity of these discrete conformal structures.
Recently, the author \cite{Xu rigidity} proved the rigidity of Glickenstein's discrete conformal structures on triangulated surfaces.
This confirms Glickenstein's conjecture \cite{G5} on the rigidity of Glickenstein's discrete conformal structures on triangulated surfaces,
and generalizes the rigidity results obtained by Bobenko-Pinkall-Springborn \cite{BPS}, Guo \cite{Guo},
Guo-Luo \cite{GL}, Luo \cite{L1, L4} and the author \cite{X1,X3}
for different special types of discrete conformal structures on triangulated surfaces.
In this paper, we study the deformation of Glickenstein's discrete conformal structures on triangulated surfaces.

Finding piecewise constant curvature metrics on polyhedral surfaces with prescribed combinatorial curvatures
is an important problem in discrete conformal geometry. This has lots of applications in theory and practice. See \cite{CL, Gu-Yau, SWGL, ZG} and others for example.
An effective approach to this problem is the combinatorial curvature flows,
including the combinatorial Ricci flow, the combinatorial Yamabe flow, the combinatorial Calabi flow and others.
The properties of these combinatorial curvature flows for Thurston's circle packings, Bowers-Stephenson's inversive distance circle packings and
Luo's vertex scalings have been extensively studied in the literature.
See \cite{BL, BL2, Ge-thesis,Ge,GH1,GJ0, GJ1, GJ3, GX3,CL,LZ, L1, ZX,GLSW,GGLSW,Wu, WX, XZ2023, XZ2023b} and others.
However, the properties of these combinatorial curvature flows for Glickenstein's discrete conformal structures are not well studied.
In \cite{ZGZLYG}, Zhang-Guo-Zeng-Luo-Yau-Gu introduced the combinatorial Ricci flow
for Glickenstein's discrete conformal structures.
Due to the singularities possibly developed along the combinatorial Ricci flow,
the longtime behavior of the combinatorial Ricci flow has been a difficult problem.
In this paper, we extend the combinatorial curvature by constants to handle the potential singularities along the combinatorial Ricci flow.
It is proved that the solution of the combinatorial Ricci flow can be uniquely extended
to exist for all time, and the solution of the extended combinatorial Ricci flow converges exponentially fast
if the discrete conformal factor with the prescribed combinatorial curvature exists.
We further introduce the combinatorial Calabi flow for Glickenstein's discrete conformal structures on triangulated surfaces,
and establish the basic properties of the combinatorial Calabi flow.
These combinatorial curvature flows provide effective algorithms for finding piecewise constant curvature metrics on polyhedral surfaces with prescribed combinatorial curvatures.

Suppose $(M, \mathcal{T})$ is a connected closed triangulated surface with a triangulation $\mathcal{T}$,
which is the quotient of a finite disjoint union of triangles by identifying all the edges of triangles
in pair by homeomorphisms. We use $V, E, F$ to denote the set of vertices, unoriented edges and faces
in the triangulation $\mathcal{T}$ respectively. For simplicity, we use one index to denote a vertex (such as $i\in V$),
two indices to denote an edge (such as $\{ij\}\in E$) and three indices to denote a triangle (such as $\{ijk\}\in F$).
We further use $f_i=f(i)$ for a function $f: V\rightarrow \mathbb{R}$, $f_{ij}=f(\{ij\})$ for a function $f: E\rightarrow \mathbb{R}$,
and $f_{ijk}=f(\{ijk\})$ for a function $f: F\rightarrow \mathbb{R}$ for simplicity.
Denote  the set of positive real numbers as $\mathbb{R}_{>0}$ and $|V|=N$.

\begin{definition}[\cite{L5}]\label{defn of polyhedral surface}
A polyhedral surface $(M, \mathcal{T}, l)$ with background geometry $\mathbb{G}$ ($\mathbb{G}=\mathbb{E}^2, \mathbb{H}^2$ or $\mathbb{S}^2$)
is a triangulated surface $(M, \mathcal{T})$ with a map $l: E\rightarrow \mathbb{R}_{>0}$ such that
any face $\{ijk\}\in F$ can be embedded in $\mathbb{G}$ as a nondegenerate triangle
with edge lengths $l_{ij}, l_{ik}, l_{jk}$ given by $l$.
The map $l: E\rightarrow \mathbb{R}_{>0}$ is called a Euclidean (hyperbolic or spherical respectively) polyhedral metric
if $\mathbb{G}=\mathbb{E}^2$ ($\mathbb{G}=\mathbb{H}^2$ or $\mathbb{G}=\mathbb{S}^2$ respectively).
\end{definition}

The nondegenerate condition for the face $\{ijk\}\in F$ in Definition \ref{defn of polyhedral surface}
is equivalent to the edge lengths $l_{ij}, l_{ik}, l_{jk}$ satisfy the triangle inequalities
($l_{ij}+l_{ik}+l_{jk}<2\pi$ additionally if $\mathbb{G}=\mathbb{S}^2$).
Intuitively, a polyhedral surface with background geometry $\mathbb{G}$ ($\mathbb{G}=\mathbb{E}^2, \mathbb{H}^2$ or $\mathbb{S}^2$) can be obtained
by gluing triangles in $\mathbb{G}$ isometrically along the edges in pair.
For polyhedral surfaces, there may exist conic singularities at the vertices,
which can be described by combinatorial curvatures.
The combinatorial curvature is a map $K: V\rightarrow (-\infty, 2\pi)$ that assigns the vertex $i\in V$
$2\pi$ less the sum of inner angles at $i$, i.e.
\begin{equation}
K_i=2\pi-\sum_{\{ijk\}\in F}\theta_{i}^{jk},
\end{equation}
where $\theta_{i}^{jk}$ is the inner angle at $i$ in the triangle $\{ijk\}\in F$.
The combinatorial curvature $K$ satisfies the following discrete Guass-Bonnet formula (\cite{CL} Proposition 3.1)
\begin{equation}\label{discrete GB formula}
\sum_{i\in V}K_i=2\pi\chi(M)-\lambda \text{Area}(M)
\end{equation}
 on a closed polyhedral surface $M$, where $\lambda=-1,0,1$ for $\mathbb{G}=\mathbb{H}^2, \mathbb{E}^2, \mathbb{S}^2$ respectively and $\text{Area}(M)$ is the area of the polyhedral surface $M$.

\begin{definition}[\cite{G3, GT}]\label{defn of discrete conformal structure}
Suppose $(M, \mathcal{T})$ is a triangulated connected closed surface and $\varepsilon: V\rightarrow \{-1, 0, 1\}$,
$\eta: E\rightarrow \mathbb{R}$ are two weights defined on the vertices and edges respectively.
A discrete conformal structure on the weighted triangulated surface $(M, \mathcal{T}, \varepsilon, \eta)$ with background geometry $\mathbb{G}$
is composed of  the maps $f: V\rightarrow \mathbb{R}$ such that
\begin{description}
  \item[(1)] the edge length $l_{ij}$ for the edge $\{ij\}\in E$ is given by
  \begin{equation}\label{defn of Euclidean length}
  \begin{aligned}
  l_{ij}=\sqrt{\varepsilon_i e^{2f_i}+\varepsilon_j e^{2f_j}+2\eta_{ij}e^{f_i+f_j}}
  \end{aligned}
  \end{equation}
  for $\mathbb{G}=\mathbb{E}^2$,
  \begin{equation}\label{defn of hyperbolic length}
  \begin{aligned}
  l_{ij}=\cosh^{-1}\left(\sqrt{(1+\varepsilon_ie^{2f_i})(1+\varepsilon_je^{2f_j})}+\eta_{ij}e^{f_i+f_j}\right)
  \end{aligned}
  \end{equation}
  for $\mathbb{G}=\mathbb{H}^2$ and
  \begin{equation}\label{defn of spherical length}
  \begin{aligned}
  l_{ij}=\cos^{-1}\left(\sqrt{(1-\varepsilon_ie^{2f_i})(1-\varepsilon_je^{2f_j})}-\eta_{ij}e^{f_i+f_j}\right)
  \end{aligned}
  \end{equation}
  for $\mathbb{G}=\mathbb{S}^2$;
  \item[(2)]  the edge length function  $l: E\rightarrow \mathbb{R}_{>0}$ defined by
  (\ref{defn of Euclidean length}), (\ref{defn of hyperbolic length}), (\ref{defn of spherical length})
  is a Euclidean, hyperbolic and spherical polyhedral metric on $(M, \mathcal{T})$ respectively.
\end{description}
The function $f: V\rightarrow \mathbb{R}$ is called a discrete conformal factor.
\end{definition}


\begin{remark}
Glickenstein's discrete conformal structures in Definition \ref{defn of discrete conformal structure} contain
different types of circle packings and Luo's vertex scalings as special cases.
Furthermore, the discrete conformal structures in Definition \ref{defn of discrete conformal structure}
contain the mixed type of discrete conformal structures. For example, it contains the type with
$\varepsilon_i=0$ for some vertices $i\in V$ and $\varepsilon_j=1$ for the other vertices $j\in V$.
In this paper, we focus on Glickenstein's Euclidean and hyperbolic discrete conformal structures
in Definition \ref{defn of discrete conformal structure}.
\end{remark}

For the discrete conformal structures in Definition \ref{defn of discrete conformal structure},
the combinatorial Ricci flow was introduced by Zhang-Guo-Zeng-Luo-Yau-Gu \cite{ZGZLYG}.
For simplicity, set
\begin{equation}\label{u Euclidean}
\begin{aligned}
u_i=f_i
\end{aligned}
\end{equation}
for any $i\in V$ in the Euclidean background geometry and
\begin{equation}\label{u hyperbolic}
\begin{aligned}
u_i=\left\{
      \begin{array}{ll}
        f_i, & \hbox{ if $\varepsilon_i=0$,} \\
        \frac{1}{2}\log \left|\frac{\sqrt{1+\varepsilon_ie^{2f_i}}-1}{\sqrt{1+\varepsilon_ie^{2f_i}}+1}\right|, & \hbox{ if $\varepsilon_i\neq 0$,}
      \end{array}
    \right.
\end{aligned}
\end{equation}
for the hyperbolic background geometry. For simplicity, $u: V\rightarrow \mathbb{R}$ is also called a discrete conformal factor in the following.

\begin{definition}[\cite{ZGZLYG}]\label{defn of combinatorial Ricci flow}
Suppose $(M, \mathcal{T}, \varepsilon, \eta)$ is a weighted triangulated connected closed surface
with weights $\varepsilon: V\rightarrow \{-1, 0, 1\}$ and $\eta: E\rightarrow \mathbb{R}$.
The combinatorial Ricci flow for Glickenstein's discrete conformal structures on polyhedral surfaces is defined as
\begin{equation}\label{CRF equ}
\begin{aligned}
\frac{du_i}{dt}=-K_i
\end{aligned}
\end{equation}
for the Euclidean and hyperbolic background geometry.
\end{definition}

The normalized combinatorial Ricci flow for Glickenstein's Euclidean discrete conformal structures  is
\begin{equation}\label{CRF equ nomalized}
\begin{aligned}
\frac{du_i}{dt}=K_{av}-K_i,
\end{aligned}
\end{equation}
where $K_{av}=\frac{2\pi \chi(M)}{N}$ is the average combinatorial curvature.

\begin{remark}
The combinatorial Ricci flow in Definition \ref{defn of combinatorial Ricci flow}
generalizes the known forms of combinatorial Ricci flow or combinatorial Yamabe flow for different types of discrete conformal
structures on polyhedral surfaces.
If $\varepsilon\equiv1$, the combinatorial Ricci flow in Definition \ref{defn of combinatorial Ricci flow}
is reduced to the combinatorial Ricci flows for different types of circle packings on polyhedral surfaces in \cite{CL,GJ1,GJ3}.
If $\varepsilon\equiv0$, the combinatorial Ricci flow in Definition \ref{defn of combinatorial Ricci flow}
is reduced to Luo's combinatorial Yamabe flow for the vertex scalings on polyhedral surfaces in \cite{L1}.
The combinatorial Ricci flow in Definition \ref{defn of combinatorial Ricci flow}
further includes the mixed type case that $\varepsilon_i=1$ for some vertices $i\in V_1\neq \emptyset$ and $\varepsilon_j=0$
for the other vertices $j\in V\setminus V_1\neq \emptyset$. The combinatorial Ricci flow in Definition \ref{defn of combinatorial Ricci flow}
can also be defined for spherical background geometry with $u_i$ satisfying
$\frac{\partial f_i}{\partial u_i}=\sqrt{1-\varepsilon_i e^{2f_i}}$. Please refer to \cite{GT,ZGZLYG} for more details.
\end{remark}

The combinatorial Ricci flows (\ref{CRF equ}) and (\ref{CRF equ nomalized}) on triangulated surfaces
may develop singularities, which correspond to
the triangles in the polyhedral surfaces degenerate.
To handle the potential singularities along the combinatorial Ricci flows (\ref{CRF equ}) and (\ref{CRF equ nomalized}),
we extend the combinatorial curvature by constants and then extend the combinatorial Ricci flow through the singularities.
We have the following result on the longtime existence and convergence for the solutions of the extended combinatorial Ricci flow
in the Euclidean and hyperbolic background geometry.

\begin{theorem}\label{convergence of extended CRF introduction}
Suppose $(M, \mathcal{T}, \varepsilon, \eta)$ is a weighted triangulated connected closed surface with
the weights $\varepsilon: V\rightarrow \{0, 1\}$ and $\eta: E\rightarrow \mathbb{R}$ satisfying
the structure conditions
\begin{equation}\label{structure condition 1}
\begin{aligned}
\varepsilon_s \varepsilon_t +\eta_{st}>0, \ \ \forall \{st\}\in E
\end{aligned}
\end{equation}
and
\begin{equation}\label{structure condition 2}
\begin{aligned}
\varepsilon_q\eta_{st}+\eta_{qs}\eta_{qt}\geq 0, \ \ \{q, s, t\}= \{i, j, k\}, \forall \{ijk\}\in F.
\end{aligned}
\end{equation}
\begin{description}
  \item[(a)] The solution of normalized combinatorial Ricci flow (\ref{CRF equ nomalized}) in the Euclidean background geometry
and the solution of the combinatorial Ricci flow (\ref{CRF equ})
in the hyperbolic background geometry can be extended
to exist for all time for any initial discrete conformal factor on $(M, \mathcal{T}, \varepsilon, \eta)$.
  \item[(b)] The solution of the extended combinatorial Ricci flow is unique for any initial
discrete conformal factor.
  \item[(c)] If there exists a Euclidean discrete conformal factor with
constant combinatorial curvature on $(M, \mathcal{T}, \varepsilon, \eta)$,
the solution of the extended normalized Euclidean combinatorial Ricci flow converges exponentially fast for
any initial Euclidean discrete conformal factor;
If there exists a hyperbolic discrete conformal factor with zero combinatorial curvature on $(M, \mathcal{T}, \varepsilon, \eta)$,
the solution of the extended hyperbolic combinatorial Ricci flow
converges exponentially fast for any initial hyperbolic discrete conformal factor.
\end{description}
\end{theorem}

\begin{remark}
If $\varepsilon\equiv1$, the result in Theorem \ref{convergence of extended CRF introduction}
is reduced to the convergence results for the combinatorial Ricci flow for Thurston's circle packings obtained in \cite{CL}
and for Bowers-Stephenson's inversive distance circle packings obtained in \cite{GJ1,GJ3}.
If $\varepsilon\equiv0$, the result in Theorem \ref{convergence of extended CRF introduction}
is reduced to the convergence result for Luo's combinatorial Yamabe flow for the vertex scalings obtained in \cite{GJ0}.
The results in Theorem \ref{convergence of extended CRF introduction}
further contain the case of mixed type that $\varepsilon_i=1$ for some vertices $i\in V_1\neq \emptyset$ and $\varepsilon_j=0$
for the other vertices $j\in V\setminus V_1\neq \emptyset$.
\end{remark}

Combinatorial Calabi flow is another effective combinatorial curvature flow for finding
piecewise constant curvature metrics on polyhedral surfaces with prescribed combinatorial curvatures.
It was first introduced by Ge \cite{Ge-thesis} (see also \cite{Ge}) for Thurston's Euclidean circle packings.
The combinatorial Calabi flow for Glickenstein's discrete conformal structures on polyhedral surfaces is defined as follows.

\begin{definition}\label{defn of combinatorial Calabi flow}
Suppose $(M, \mathcal{T}, \varepsilon, \eta)$ is a weighted triangulated connected closed surface
with weights $\varepsilon: V\rightarrow \{-1, 0, 1\}$ and $\eta: E\rightarrow \mathbb{R}$.
The combinatorial Calabi flow for Glickenstein's discrete conformal structures on polyhedral surfaces is defined as
\begin{equation}\label{CCF equ}
\begin{aligned}
\frac{du_i}{dt}=\Delta K_i
\end{aligned}
\end{equation}
for the Euclidean and hyperbolic background geometry, where $\Delta=-\frac{\partial (K_1, \cdots, K_N)}{\partial (u_1, \cdots, u_N)}$ is the combinatorial Laplace operator.
\end{definition}

\begin{remark}
The combinatorial Calabi flow introduced in Definition \ref{defn of combinatorial Calabi flow}
unifies the known forms of combinatorial Calabi flows for different special types of discrete conformal structures on polyhedral surfaces.
If $\varepsilon\equiv1$, the combinatorial Calabi flow (\ref{CCF equ})
is reduced to the combinatorial Calabi flow for circle packings on polyhedral surfaces introduced in \cite{Ge-thesis,Ge,GX3}.
If $\varepsilon\equiv0$, the combinatorial Calabi flow (\ref{CCF equ})
is reduced to the combinatorial Calabi flow for Luo's vertex scalings on polyhedral surfaces introduced in \cite{Ge-thesis}.
The combinatorial Calabi flow in Definition \ref{defn of combinatorial Calabi flow}
further contains the mixed type case that $\varepsilon_i=1$ for some vertices $i\in V_1\neq \emptyset$ and $\varepsilon_j=0$
for the other vertices $j\in V\setminus V_1\neq \emptyset$.
Similar to the combinatorial Ricci flow, the combinatorial Calabi flow (\ref{CCF equ})
can also be defined for
spherical background geometry with $u_i$ satisfying
$\frac{\partial f_i}{\partial u_i}=\sqrt{1-\varepsilon_i e^{2f_i}}$.
\end{remark}

The combinatorial Calabi flow is a negative gradient flow of the combinatorial Calabi energy $\mathcal{C}=\frac{1}{2}\sum_{i=1}^NK_i^2$.
This fact was first observed by Ge \cite{Ge-thesis,Ge} in the case of Thurston's circle packings.
We have the following result on the longtime behavior of the combinatorial Calabi flow (\ref{CCF equ}).

\begin{theorem}\label{convergence of CCF introduction}
Suppose $(M, \mathcal{T}, \varepsilon, \eta)$ is a weighted triangulated connected closed surface with
the weights $\varepsilon: V\rightarrow \{0, 1\}$ and $\eta: E\rightarrow \mathbb{R}$ satisfying
the structure conditions (\ref{structure condition 1}) and (\ref{structure condition 2}).
\begin{description}
  \item[(a)] If the solution of the combinatorial Calabi flow (\ref{CCF equ}) converges to a nondegenerate discrete conformal
factor $u^*$, then $u^*$ has constant combinatorial curvature.
  \item[(b)] If there exists a discrete conformal factor $u^*$ with constant combinatorial curvature ($\frac{2\pi\chi(M)}{N}$ for the Euclidean background geometry and $0$ for the hyperbolic background geometry), then there exists a real number $\delta>0$
such that if the initial value $u(0)$ of the combinatorial Calabi flow (\ref{CCF equ}) satisfies $||u(0)-u^*||<\delta$,
the solution of the combinatorial Calabi flow (\ref{CCF equ}) exists for all time and converges exponentially fast to $u^*$.
\end{description}
\end{theorem}

The paper is organized as follows.
In Section \ref{section 2}, we study the Euclidean combinatorial Ricci flow and Euclidean combinatorial Calabi flow on polyhedral surfaces and
prove a generalization of Theorem \ref{convergence of extended CRF introduction} and Theorem \ref{convergence of CCF introduction} in the Euclidean background geometry.
In Section \ref{section 3}, we study the hyperbolic combinatorial Ricci flow and hyperbolic combinatorial Calabi flow on polyhedral surfaces and
prove a generalization of Theorem \ref{convergence of extended CRF introduction} and Theorem \ref{convergence of CCF introduction} in the hyperbolic background geometry.
In Section \ref{section 4}, we discuss some open problems.
\\
\\
\textbf{Acknowledgements}\\[8pt]
The author thanks the referee for reading the paper carefully and pointing out a simplification of the proof of Theorem \ref{Euclidean uniqueness}.
The research of the author is partially supported by Fundamental Research Funds for the Central Universities under Grant no. 2042020kf0199.\\

\section{Euclidean combinatorial curvature flows}\label{section 2}
For the further applications, we study the following modified Euclidean combinatorial Ricci flow
\begin{equation}\label{modified CRF equ}
\begin{aligned}
\frac{du_i}{dt}=\overline{K}_i-K_i
\end{aligned}
\end{equation}
and
modified Euclidean combinatorial Calabi flow
\begin{equation}\label{modified CCF equ}
\begin{aligned}
\frac{du_i}{dt}=\Delta (K-\overline{K})_i
\end{aligned}
\end{equation}
for Glickenstein's Euclidean discrete conformal structures on triangulated surfaces,
where $\overline{K}: V\rightarrow (-\infty, 2\pi)$ is a function defined on the vertices
with $\sum_{i=1}^N\overline{K}_i=2\pi\chi(M)$.
The modified combinatorial Ricci flow (\ref{modified CRF equ})
and modified combinatorial Calabi flow (\ref{modified CCF equ})
are generalizations of the normalized combinatorial Ricci flow (\ref{CRF equ nomalized})
and the combinatorial Calabi flow (\ref{CCF equ}) respectively,
and can be used to study the prescribed combinatorial curvature problem for Glickenstein's Euclidean
discrete conformal structures on polyhedral surfaces.

We first recall some basic properties of Glickenstein's Euclidean discrete conformal structures on triangulated surfaces.
For more details on these properties, please refer to \cite{GT,Xu rigidity}.
Suppose $\sigma=\{ijk\}$ is a topological triangle, we use $V_\sigma=\{i,j,k\}$ to denote its vertices and use
$E_\sigma=\{\{ij\}, \{ik\}, \{jk\}\}$ to denote its edges in the following.

\begin{theorem}[\cite{Xu rigidity}]\label{theorem admissible space Euclidean}
Suppose $\sigma=\{ijk\}$ is a topological triangle with
two weights $\varepsilon: V_\sigma\rightarrow \{0, 1\}$ and $\eta: E_\sigma\rightarrow \mathbb{R}$
satisfying the structure conditions (\ref{structure condition 1}) and (\ref{structure condition 2}).
Then the admissible space $\Omega_{ijk}^E(\eta)$ of
nondegerate Euclidean discrete conformal factors $u$ for the triangle $\sigma=\{ijk\}$
is a nonempty simply connected open set whose boundary components are analytic.
Furthermore, the admissible space $\Omega_{ijk}^E(\eta)$ in the parameter $u$ can be written as
\begin{equation*}
\begin{aligned}
\Omega_{ijk}^E(\eta)=\mathbb{R}^3\setminus\sqcup_{\alpha\in \Lambda}V_\alpha,
\end{aligned}
\end{equation*}
where
$\Lambda=\{q\in \{i,j,k\}|A_q=\eta_{st}^2-\varepsilon_s\varepsilon_t>0, \{q, s, t\}=\{i, j, k\}\}$,
$\sqcup_{\alpha\in \Lambda}V_\alpha$ is a disjoint union of $V_\alpha$ and
$V_\alpha$ is a closed region in $\mathbb{R}^3$ bounded by an analytical function defined on $\mathbb{R}^2$.
\end{theorem}

\begin{remark}\label{derivative tends infty Euclidean}
If $(\eta_{ij}, \eta_{ik}, \eta_{jk})$ defined on the triangle $\{ijk\}$
satisfies the structure conditions (\ref{structure condition 1}) and (\ref{structure condition 2}) and
 $(u_i, u_j, u_k)\in \Omega_{ijk}^E(\eta)$ tends to a point $(\overline{u}_i, \overline{u}_j, \overline{u}_k)\in \partial V_i$ with $V_i\neq \emptyset$,
 we have $\frac{\partial \theta_i}{\partial u_j}\rightarrow +\infty$, $\frac{\partial \theta_i}{\partial u_k}\rightarrow +\infty$ and
 $\frac{\partial \theta_i}{\partial u_i}\rightarrow -\infty$.
\end{remark}

\begin{theorem}\label{Euclidean Jacobian negativity}
Suppose $\sigma=\{ijk\}$ is a topological triangle with
two weights $\varepsilon: V_\sigma\rightarrow \{0, 1\}$ and $\eta: E_\sigma\rightarrow \mathbb{R}$
satisfying the structure conditions (\ref{structure condition 1}) and (\ref{structure condition 2}).
Then the Jacobian matrix $\Lambda^E_{ijk}=\frac{\partial (\theta_i, \theta_j, \theta_k)}{\partial ( u_i, u_j, u_k)}$
is symmetric, negative semi-definite with rank $2$ and  has kernel $\{t(1,1,1)^T|t\in \mathbb{R}\}$ for any nondegenerate
Euclidean discrete conformal factor on the triangle $\{ijk\}$.
\end{theorem}

\begin{corollary}\label{Euclidean curvature Jacobian positivity}
Suppose $(M, \mathcal{T}, \varepsilon, \eta)$ is a weighted triangulated surface with
the weights $\varepsilon: V\rightarrow \{0, 1\}$ and $\eta: E\rightarrow \mathbb{R}$
satisfying the structure conditions (\ref{structure condition 1}) and (\ref{structure condition 2}).
Then the Jacobian matrix $\Lambda^E=\frac{\partial (K_1,\cdots,K_N)}{\partial (u_1,\cdots, u_N)}$
is symmetric and positive semi-definite with rank $N-1$ and has kernel $\{t\mathbf{1}\in \mathbb{R}^N|t\in \mathbb{R}\}$
for all nondegenerate Euclidean discrete conformal factors on $(M, \mathcal{T}, \varepsilon, \eta)$.
\end{corollary}

Theorem \ref{theorem admissible space Euclidean} and Theorem \ref{Euclidean Jacobian negativity} imply the following function
\begin{equation}\label{Euclidean Ricci energy function for triangle}
\begin{aligned}
\mathcal{E}_{ijk}(u_i,u_j,u_k)=\int_{(\overline{u}_i, \overline{u}_j, \overline{u}_k)}^{(u_i, u_j, u_k)}\theta_idu_i+\theta_jdu_j+\theta_kdu_k
\end{aligned}
\end{equation}
is a well-defined smooth locally concave function on $\Omega^E_{ijk}(\eta)$ with $\nabla_{u_i}\mathcal{E}_{ijk}=\theta_i$ and
$\mathcal{E}_{ijk}(u_i+t,u_j+t,u_k+t)=\mathcal{E}_{ijk}(u_i,u_j,u_k)+t\pi$.
The function $\mathcal{E}_{ijk}$ is called the Ricci energy function for the triangle $\{ijk\}$.
Set
\begin{equation}\label{Euclidean Ricci energy function}
\begin{aligned}
\mathcal{E}(u_1, \cdots, u_N)=2\pi\sum_{i\in V}u_i-\sum_{\{ijk\}\in F}\mathcal{E}_{ijk}(u_i,u_j,u_k)
\end{aligned}
\end{equation}
to be the Ricci energy function defined on the admissible space $\Omega^E$ of nondegenerate Euclidean discrete conformal factors
on $(M, \mathcal{T}, \varepsilon, \eta)$.
Then $\mathcal{E}$ is a locally convex function defined on $\Omega^E$ with
$\mathcal{E}(u_1+t, \cdots, u_N+t)=\mathcal{E}(u_1, \cdots, u_N)+2t\pi\chi(M)$ and $\nabla_{u_i} \mathcal{E}=K_i$
by Corollary \ref{Euclidean curvature Jacobian positivity}.

\begin{lemma}\label{CRF and CCF are gradient flows}
Suppose $(M, \mathcal{T}, \varepsilon, \eta)$ is a weighted triangulated connected closed surface with
the weights $\varepsilon: V\rightarrow \{0, 1\}$ and $\eta: E\rightarrow \mathbb{R}$ satisfying
the structure conditions (\ref{structure condition 1}) and (\ref{structure condition 2}).
The modified Euclidean combinatorial Ricci flow (\ref{modified CRF equ})
and modified Euclidean combinatorial Calabi flow (\ref{modified CCF equ}) on $(M, \mathcal{T}, \varepsilon, \eta)$
are negative gradient flows.
\end{lemma}
\proof
Set
$\mathcal{H}(u)=\mathcal{E}(u_1, \cdots, u_N)-\sum_{i=1}^N\overline{K}_iu_i,$
where $\mathcal{E}(u_1, \cdots, u_N)$ is the Ricci energy function defined by (\ref{Euclidean Ricci energy function}).
Then $\nabla_{u_i}\mathcal{H}=K_i-\overline{K}_i$, which implies the modified Euclidean combinatorial Ricci flow (\ref{modified CRF equ})
is a negative gradient flow of $\mathcal{H}(u)$.

Set $\mathcal{C}(u)=\frac{1}{2}||\overline{K}-K||^2=\frac{1}{2}\sum_{i=1}^N(\overline{K}_i-K_i)^2$.
By direct calculations, we have $\nabla_{u_i}\mathcal{C}=-\Delta (K-\overline{K})_i$,
which implies the modified Euclidean combinatorial Calabi flow (\ref{modified CCF equ}) is a negative gradient flow of $\mathcal{C}(u)$.
\qed

As the modified Euclidean combinatorial Ricci flow (\ref{modified CRF equ})
and modified Euclidean combinatorial Calabi flow (\ref{modified CCF equ}) are ODE systems,
the short time existence of the solutions are ensured by the standard ODE theory.
We further have the following result on the longtime existence and convergence for the solutions of
the modified combinatorial Ricci flow (\ref{modified CRF equ})
and the modified combinatorial Calabi flow (\ref{modified CCF equ})
for initial values with small energy,
which is a slight generalization of Theorem \ref{convergence of CCF introduction} in the Euclidean background geometry.

\begin{theorem}\label{stability Euclidean}
Suppose $(M, \mathcal{T}, \varepsilon, \eta)$ is a weighted triangulated connected closed surface with
the weights $\varepsilon: V\rightarrow \{0, 1\}$ and $\eta: E\rightarrow \mathbb{R}$ satisfying
the structure conditions (\ref{structure condition 1}) and (\ref{structure condition 2}).
\begin{description}
  \item[(a)] If the solution of the modified Euclidean combinatorial Ricci flow (\ref{modified CRF equ})
or modified Euclidean  combinatorial Calabi flow (\ref{modified CCF equ}) converges to a nondegenerate discrete conformal
factor $\overline{u}$, then the combinatorial curvature of the Euclidean polyhedral metric
determined by the discrete conformal factor $\overline{u}$ is $\overline{K}$.
  \item[(b)] If there exists a nondegenerate Euclidean discrete conformal factor $\overline{u}$ with combinatorial curvature $\overline{K}$,
  then there exists a real number $\delta>0$ such that if the initial value $u(0)$ of
   the modified Euclidean combinatorial Ricci flow (\ref{modified CRF equ}) (the modified Euclidean combinatorial Calabi flow (\ref{modified CCF equ}) respectively)
   satisfies $||u(0)-\overline{u}||<\delta$, the solution of the modified Euclidean combinatorial Ricci flow (\ref{modified CRF equ})
  (the modified Euclidean combinatorial Calabi flow (\ref{modified CCF equ}) respectively) exists for all time and converges exponentially fast to $\overline{u}$.
\end{description}
\end{theorem}
\proof
The proof for part (a) is direct. We just prove part (b).
For the modified Euclidean combinatorial Ricci flow (\ref{modified CRF equ}),
by direct calculations, we have
\begin{equation}\label{invariant of sum of ui along CRF}
\frac{d(\sum_{i=1}^N u_i)}{dt}=\sum_{i=1}^N(\overline{K}_i-K_i)=2\pi\chi(M)-2\pi\chi(M)=0,
\end{equation}
where the discrete Gauss-Bonnet formula (\ref{discrete GB formula}) for the Euclidean background geometry and the assumption $\sum_{i=1}^N\overline{K}_i=2\pi\chi(M)$
are used in the second equality.
This implies $\sum_{i=1}^N u_i$ is invariant along the modified combinatorial Ricci flow (\ref{modified CRF equ}).
Without loss of generality, assume $\sum_{i=1}^Nu_i(0)=0$.  Set
$$\Sigma=\{u\in \mathbb{R}^N|\sum_{i=1}^Nu_i=0\}.$$
By (\ref{invariant of sum of ui along CRF}),
the solution $u(t)$ of the modified combinatorial Ricci flow (\ref{modified CRF equ}) stays in the hyperplane $\Sigma$.
Set $\Gamma(u)=\overline{K}-K$ for the modified combinatorial Ricci flow (\ref{modified CRF equ}).
Then $\overline{u}$ is an equilibrium point of the system (\ref{modified CRF equ}) and
$D\Gamma(\overline{u})=-\frac{\partial (K_1,\cdots,K_N)}{\partial (u_1,\cdots, u_N)}$
is negative semi-definite with kernel space $\{t\mathbf{1}\in \mathbb{R}^N|t\in \mathbb{R}\}$ by Corollary \ref{Euclidean curvature Jacobian positivity}.
Note that the solution $u(t)$ of the modified combinatorial Ricci flow (\ref{modified CRF equ}) stays in the hyperplane $\Sigma$,
the normal vector of which generates the kernel space $\{t\mathbf{1}\in \mathbb{R}^N|t\in \mathbb{R}\}$ of $D\Gamma(\overline{u})$.
This implies $\overline{u}$ is a local attractor of the system (\ref{modified CRF equ}).
Then the longtime existence and exponential convergence of the solution of (\ref{modified CRF equ})
follows from the Lyapunov Stability Theorem (\cite{P}, Chapter 5).

For the modified Euclidean combinatorial Calabi flow (\ref{modified CCF equ}), we have
\begin{equation*}
\frac{d(\sum_{i=1}^N u_i)}{dt}=\sum_{i=1}^N\Delta(\overline{K}-K)_i
=\sum_{j=1}^N\sum_{i=1}^N\Lambda^E_{ij}(\overline{K}-K)_j=0
\end{equation*}
by the kernel space of $\Lambda^E$ is $\{t(1,\cdots, 1)\in \mathbb{R}^N|t\in \mathbb{R}\}$ in Corollary \ref{Euclidean curvature Jacobian positivity}.
This implies $\sum_{i=1}^N u_i$ is invariant along the flow (\ref{modified CCF equ}).
Set $\Gamma(u)=\Delta (K-\overline{K})$. Then $\overline{u}$ is an equilibrium point of the system (\ref{modified CCF equ}) and
$D\Gamma(\overline{u})=-\left(\frac{\partial (K_1,\cdots,K_N)}{\partial (u_1,\cdots, u_N)}\right)^2$
is negative semi-definite with kernel space $\{t(1,\cdots, 1)\in \mathbb{R}^N|t\in \mathbb{R}\}$ by Corollary \ref{Euclidean curvature Jacobian positivity}.
The rest of the proof is the same as that for the modified Euclidean combinatorial Ricci flow (\ref{modified CRF equ}), we omit the details here.
\qed

For general initial value, the modified Euclidean  combinatorial Ricci flow (\ref{modified CRF equ})
and the modified Euclidean  combinatorial Calabi flow (\ref{modified CCF equ}) may develop singularities,
which correspond to the triangles in the triangulation degenerate or the discrete conformal factor $f$
tends to infinity along the combinatorial curvature flows.
For the modified Euclidean  combinatorial Ricci flow (\ref{modified CRF equ}), we can extend it through the singularities to ensure the
longtime existence and convergence of the solution for general initial value.
We first recall the following extension introduced by the author \cite{Xu rigidity} for Glickenstein's Euclidean discrete conformal structures.

\begin{lemma}\label{Euclidean extension}
Suppose $\sigma=\{ijk\}$ is a topological triangle with
two weights $\varepsilon: V_\sigma\rightarrow \{0, 1\}$ and $\eta: E_\sigma\rightarrow \mathbb{R}$
satisfying the structure conditions (\ref{structure condition 1}) and (\ref{structure condition 2}).
Then the inner angle functions $\theta_i, \theta_j, \theta_k$ defined on $\Omega^E_{ijk}(\eta)$
can be extended to be continuous functions $\widetilde{\theta}_i, \widetilde{\theta}_j, \widetilde{\theta}_k$
defined on $\mathbb{R}^3$ by setting
\begin{equation}\label{extension of theta_i}
\begin{aligned}
\widetilde{\theta}_i(u_i, u_j, u_k)=\left\{
                                      \begin{array}{ll}
                                        \theta_i, & \hbox{if $(u_i, u_j, u_k)\in \Omega^E_{ijk}(\eta)$;} \\
                                        \pi, & \hbox{if $(u_i, u_j, u_k)\in V_i$;} \\
                                        0, & \hbox{otherwise.}
                                      \end{array}
                                    \right.
\end{aligned}
\end{equation}
\end{lemma}

By the extension in Lemma \ref{Euclidean extension}, we can extend the combinatorial curvature $K$ defined on $\Omega^E$ to be the generalized combinatorial curvature $\widetilde{K}_i:=2\pi-\sum\widetilde{\theta}_i$ defined for any $u\in \mathbb{R}^N$.
We call $u\in \mathbb{R}^N$ as a generalized discrete conformal factor.

Recall the following definition of closed continuous $1$-forms and
extension of locally convex functions introduced Luo \cite{L4}, which is a generalization
of Bobenko-Pinkall-Spingborn's extension introduced in \cite{BPS}.
\begin{definition}[\cite{L4}, Definition 2.3]
A differential 1-form $w=\sum_{i=1}^n a_i(x)dx^i$ defined on an open set $U\subset \mathbb{R}^n$ is said to be continuous if each $a_i(x)$ is continuous on $U$. A continuous differential 1-form $w$ is called closed if $\int_{\partial \tau}w=0$ for each
triangle $\tau\subset U$.
\end{definition}

\begin{theorem}[\cite{L4}, Corollary 2.6]\label{Luo's convex extention}
Suppose $X\subset \mathbb{R}^n$ is an open convex set and $A\subset X$ is an open subset of $X$ bounded by a real analytic codimension-1 submanifold in $X$. If $w=\sum_{i=1}^na_i(x)dx_i$ is a continuous closed 1-form on $A$ so that $F(x)=\int_a^x w$ is locally convex on $A$ and each $a_i$ can be extended continuous to $X$ by constant functions to a function $\widetilde{a}_i$ on $X$, then  $\widetilde{F}(x)=\int_a^x\sum_{i=1}^n\widetilde{a}_i(x)dx_i$ is a $C^1$-smooth
convex function on $X$ extending $F$.
\end{theorem}

By Lemma \ref{Euclidean extension} and Theorem \ref{Luo's convex extention}, the locally concave function
$\mathcal{E}_{ijk}$ defined by (\ref{Euclidean Ricci energy function for triangle})
for nondegenerate Euclidean discrete conformal factors $u\in\Omega^E_{ijk}(\eta)$ on a triangle $\{ijk\}$
can be extended to be a $C^1$ smooth concave function
\begin{equation}\label{extension of Euclidean Ricci energy function for triangle}
\begin{aligned}
\widetilde{\mathcal{E}}_{ijk}(u_i,u_j,u_k)=\int_{(\overline{u}_i, \overline{u}_j, \overline{u}_k)}^{(u_i, u_j, u_k)}\widetilde{\theta}_idu_i+\widetilde{\theta}_jdu_j+\widetilde{\theta}_kdu_k
\end{aligned}
\end{equation}
defined on $\mathbb{R}^3$ with $\nabla_{u_i}\widetilde{\mathcal{E}}_{ijk}=\widetilde{\theta}_i$.
As a result, the locally convex function $\mathcal{E}$ defined by (\ref{Euclidean Ricci energy function}) for nondegenerate Euclidean
discrete conformal factors on $(M, \mathcal{T}, \varepsilon, \eta)$ can be extended to
be a $C^1$ smooth convex function
\begin{equation}\label{extended Euclidean energy function}
\begin{aligned}
\widetilde{\mathcal{E}}(u_1, \cdots, u_N)=2\pi\sum_{i\in V}u_i-\sum_{\{ijk\}\in F}\widetilde{\mathcal{E}}_{ijk}(u_i,u_j,u_k)
\end{aligned}
\end{equation}
defined on $\mathbb{R}^N$ with $\nabla_{u_i} \widetilde{\mathcal{E}}=\widetilde{K}_i=2\pi-\sum\widetilde{\theta}_i$.
Using the function $\widetilde{\mathcal{E}}$, the author \cite{Xu rigidity} proved the following rigidity of Glickenstein's Euclidean discrete conformal structures.

\begin{theorem}\label{Thm rigidity of EDCS context}
Suppose $(M, \mathcal{T}, \varepsilon, \eta)$ is a weighted triangulated surface with
the weights $\varepsilon: V\rightarrow \{0, 1\}$ and $\eta: E\rightarrow \mathbb{R}$
satisfying the structure conditions (\ref{structure condition 1}) and (\ref{structure condition 2}).
If there exists a nondegenrate Euclidean discrete conformal factor $u_A\in \Omega^E$ and
a generalized Euclidean discrete conformal factor $u_B\in \mathbb{R}^N$ such that
$K(u_A)=\widetilde{K}(u_B)$. Then $u_A=u_B+c(1,\cdots,1)$ for some constant $c\in \mathbb{R}$.
\end{theorem}

Using the extension of inner angles in Lemma \ref{Euclidean extension}, we can extend the modified Euclidean combinatorial Ricci flow (\ref{modified CRF equ}) to the following form.

\begin{definition}
Suppose $(M, \mathcal{T}, \varepsilon, \eta)$ is a weighted triangulated connected closed surface with
the weights $\varepsilon: V\rightarrow \{0, 1\}$ and $\eta: E\rightarrow \mathbb{R}$ satisfying
the structure conditions (\ref{structure condition 1}) and (\ref{structure condition 2}).
The extended modified Euclidean combinatorial Ricci flow is defined to be
\begin{equation}\label{extended modified CRF equ}
\begin{aligned}
\frac{du_i}{dt}=\overline{K}_i-\widetilde{K}_i,
\end{aligned}
\end{equation}
where $\widetilde{K}_i=2\pi-\sum_{\{ijk\}\in F}\widetilde{\theta}_i$ is an extension of the combinatorial curvature $K_i$ with
$\widetilde{\theta}_i$ given by  (\ref{extension of theta_i}).
\end{definition}

Note that the extended combinatorial curvature $\widetilde{K}$ is only a continuous function of
the generalized discrete conformal factors $u\in \mathbb{R}^N$ and does not have continuous derivatives.
Furthermore, Remark \ref{derivative tends infty Euclidean}  implies that $\widetilde{K}$ is not Lipschitz.
For such ODE systems, there may exist more than one solution.
However, we can prove the following uniqueness for the solutions of the
extended modified Euclidean combinatorial Ricci flow (\ref{extended modified CRF equ}) with any
generalized discrete conformal factor as initial value,
which is a generalization of Theorem \ref{convergence of extended CRF introduction} (b) in the Euclidean background geometry.

\begin{theorem}\label{Euclidean uniqueness}
Suppose $(M, \mathcal{T}, \varepsilon, \eta)$ is a weighted triangulated connected closed surface with
the weights $\varepsilon: V\rightarrow \{0, 1\}$ and $\eta: E\rightarrow \mathbb{R}$ satisfying
the structure conditions (\ref{structure condition 1}) and (\ref{structure condition 2}).
The solution of the extended modified Euclidean combinatorial Ricci flow (\ref{extended modified CRF equ}) is unique for any initial generalized
discrete conformal factor on $(M, \mathcal{T}, \varepsilon, \eta)$.
\end{theorem}
\proof
We prove the result by contradiction.
Suppose $u_A$ and $u_B$ are two different generalized discrete conformal factors.
Set $\widetilde{f}(t)=\widetilde{\mathcal{E}}(tu_A+(1-t)u_B)$.
By the fact that $\widetilde{\mathcal{E}}(u)$ is a $C^1$ smooth convex function defined on $\mathbb{R}^N$,
we have $\widetilde{f}(t)$ is a $C^1$ smooth convex function of $t\in [0,1]$
with $\widetilde{f}'(t)=\nabla\widetilde{\mathcal{E}}(tu_A+(1-t)u_B)\cdot (u_{A}-u_B)=\widetilde{K}(tu_A+(1-t)u_B)\cdot (u_{A}-u_B)$.
Furthermore, the derivative $\widetilde{f}'(t)$ is nondecreasing for $t\in [0,1]$.
As a result, we have
\begin{equation}\label{convex inequ}
\begin{aligned}
(\widetilde{K}(u_A)-\widetilde{K}(u_B))\cdot (u_{A}-u_B)=\widetilde{f}'(1)-\widetilde{f}'(0)\geq 0.
\end{aligned}
\end{equation}

Suppose $u_A(t)$ and $u_B(t)$, $t\in [0, T)$,  are two solutions of the extended combinatorial Ricci flow (\ref{extended modified CRF equ})
with $u_A(0)=u_B(0)$. Set $f(t)=||u_A(t)-u_B(t)||^2$. Then $f(0)=0$, $f(t)\geq 0$
and
\begin{equation*}
\begin{aligned}
\frac{df(t)}{dt}
=&2\left(\frac{du_A(t)}{dt}-\frac{du_B(t)}{dt}\right)\cdot (u_{A}(t)-u_B(t))\\
=&-(\widetilde{K}(u_A(t))-\widetilde{K}(u_B(t)))\cdot (u_{A}(t)-u_B(t))\\
\leq& 0,
\end{aligned}
\end{equation*}
where the last inequality comes from (\ref{convex inequ}). This implies $f(t)\equiv 0$ and hence $u_A(t)=u_{B}(t)$ for any $t\in [0, T)$.
\qed

For the longtime existence and convergence of the solutions of
the extended Euclidean combinatorial Ricci flow (\ref{extended modified CRF equ}),
we have the following result, which is a generalization of
Theorem \ref{convergence of extended CRF introduction} (a) (c) in the Euclidean background geometry.

\begin{theorem}\label{longtime existence and convgence of CRF Euclidean}
Suppose $(M, \mathcal{T}, \varepsilon, \eta)$ is a weighted triangulated connected closed surface with
the weights $\varepsilon: V\rightarrow \{0, 1\}$ and $\eta: E\rightarrow \mathbb{R}$ satisfying
the structure conditions (\ref{structure condition 1}) and (\ref{structure condition 2}).
The solution of the extended modified Euclidean combinatorial Ricci flow (\ref{extended modified CRF equ}) exists for all time
for any initial generalized discrete conformal factor $u(0)\in\mathbb{R}^N$ on $(M, \mathcal{T}, \varepsilon, \eta)$.
Furthermore, if there exists a nondegenerate Euclidean discrete conformal factor $\overline{u}\in\Omega^E$ with combinatorial curvature $\overline{K}$,
then the solution of the extended modified Euclidean combinatorial Ricci flow (\ref{extended modified CRF equ})
converges exponentially fast to $\overline{u}$ for any initial generalized Euclidean discrete conformal factor $u(0)\in \mathbb{R}^N$ with
$\sum_{i=1}^N u(0)=\sum_{i=1}^N \overline{u}_i$.
\end{theorem}
\proof
Suppose $u(t)$ is the solution of the extended modified Euclidean combinatorial Ricci flow (\ref{extended modified CRF equ})
with initial value $u(0)\in \mathbb{R}^N$,
then
$|\frac{du_i}{dt}|=|\overline{K}_i-\widetilde{K}_i|\leq |\overline{K}_i|+(d_i+2)\pi,$
where $d_i$ is the number of vertices adjacent to the vertex $i\in V$.
This implies $|u_i(t)|\leq |u_i(0)|+[|\overline{K}_i|+(d_i+2)\pi]t<+\infty$
for any $t\in [0, +\infty)$.
Therefore, the solution of the extended modified Euclidean combinatorial Ricci flow (\ref{extended modified CRF equ}) exists for all time.

Note that the extended inner angles $\widetilde{\theta}_i, \widetilde{\theta}_j, \widetilde{\theta}_j$ for a triangle $\{ijk\}$
in Lemma \ref{Euclidean extension} satisfy $\widetilde{\theta}_i+\widetilde{\theta}_j+\widetilde{\theta}_j=\pi$.
This implies the extended combinatorial curvature $\widetilde{K}$ satisfies the discrete Gauss-Bonnet formula
$\sum_{i=1}^N\widetilde{K}_i=2\pi\chi(M)$.
This further implies
$$\frac{d(\sum_{i=1}^N u_i)}{dt}=\sum_{i=1}^N(\overline{K}_i-\widetilde{K}_i)=0$$
along the extended modified Euclidean combinatorial Ricci flow (\ref{extended modified CRF equ}).
Therefore, $\sum_{i=1}^N u_i$ is invariant along the flow (\ref{extended modified CRF equ}).
Without loss of generality, assume $\sum_{i=1}^N u_i(0)=0$. Then the solution $u(t)$
of the extended Euclidean combinatorial Ricci flow (\ref{extended modified CRF equ})
stays in the hyperplane $\Sigma:=\{u\in \mathbb{R}^N|\sum_{i=1}^N u_i=0\}$.

Set $\widetilde{\mathcal{H}}(u)=\widetilde{\mathcal{E}}(u)-\int_{\overline{u}}^u\sum_{i=1}^N\overline{K}_idu_i$, where $\widetilde{\mathcal{E}}(u)$
is the extended Ricci energy function defined by (\ref{extended Euclidean energy function}).
Then $\widetilde{\mathcal{H}}(u)$ is a $C^1$ smooth convex function defined on $\mathbb{R}^N$ with
$\widetilde{\mathcal{H}}(u)\geq \widetilde{\mathcal{H}}(\overline{u})=0$
and $\nabla \widetilde{\mathcal{H}}(\overline{u})=K(\overline{u})-\overline{K}=0$
by the assumption $K(\overline{u})=\overline{K}$.
This further implies $\lim_{u\in \Sigma, u\rightarrow \infty}\widetilde{\mathcal{H}}(u)=+\infty$ by Corollary \ref{Euclidean curvature Jacobian positivity}
and the following property of convex functions, a proof of which can be found in \cite{GX5} (Lemma 4.6).

\begin{lemma}\label{convex func tends infty at infty}
Suppose $f(x)$ is a $C^1$ smooth convex function on $\mathbb{R}^n$ with $\nabla f(x_0)=0$ for some $x_0\in \mathbb{R}^n$,
$f(x)$ is $C^2$ smooth and strictly convex
in a neighborhood of $x_0$, then $\lim_{x\rightarrow \infty}f(x)=+\infty$.
\end{lemma}

By direct calculations, we have
\begin{equation}\label{derivative of H negative Euclidean}
\begin{aligned}
\frac{d}{dt}\widetilde{\mathcal{H}}(u(t))=\sum_{i=1}^N\nabla_{u_i} \widetilde{\mathcal{H}}\cdot \frac{du_i}{dt}=-\sum_{i=1}^N(\widetilde{K}_i-\overline{K}_i)^2\leq 0.
\end{aligned}
\end{equation}
This implies $0=\widetilde{\mathcal{H}}(\overline{u})\leq \widetilde{\mathcal{H}}(u(t))\leq \widetilde{\mathcal{H}}(u(0))$
along the extended modified Euclidean combinatorial Ricci flow (\ref{extended modified CRF equ}).
This further implies the solution $u(t)$ of the extended modified Euclidean combinatorial Ricci flow (\ref{extended modified CRF equ})
stays in a compact subset $U$ of $\Sigma$ by $\lim_{u\in \Sigma, u\rightarrow \infty}\widetilde{\mathcal{H}}(u)=+\infty$.
Therefore, $\widetilde{\mathcal{H}}(u(t))$ is bounded along the flow (\ref{extended modified CRF equ}) and the limit
$\lim_{t\rightarrow +\infty}\widetilde{\mathcal{H}}(u(t))$ exists by (\ref{derivative of H negative Euclidean}).
Taking $t_n=n$, then there exists $\xi_n\in (n, n+1)$ such that
\begin{equation}\label{K xi converges to K bar}
\begin{aligned}
\widetilde{\mathcal{H}}(u(t_{n+1}))-\widetilde{\mathcal{H}}(u(t_{n}))=-\sum_{i=1}^N\left(\widetilde{K}_i(u(\xi_n))-\overline{K}_i\right)^2\rightarrow 0,\ \text{as} \ n\rightarrow \infty.
\end{aligned}
\end{equation}
Note that $u(\xi_n)\in U\subset\subset \Sigma$, there exists a subsequence of $\{\xi_n\}$, still denoted by $\{\xi_n\}$ for simplicity,
such that $u(\xi_n)\rightarrow u^*$ for some $u^*\in U\subset\subset \Sigma$. Then $\widetilde{K}(u^*)=\overline{K}=K(\overline{u})$
by the continuity of $\widetilde{K}$ and (\ref{K xi converges to K bar}).
Therefore, $u^*=\overline{u}$ by Theorem \ref{Thm rigidity of EDCS context} and
there is a sequence $\xi_n\in (0, +\infty)$ such that $u(\xi_n)\rightarrow \overline{u}$ as $n\rightarrow \infty$.

Set $\Gamma(u)=\overline{K}-\widetilde{K}$ for the extended modified Euclidean combinatorial Ricci flow (\ref{extended modified CRF equ}).
Then $\overline{u}$ is an equilibrium point of the system (\ref{extended modified CRF equ}) and
$D\Gamma|_{\overline{u}}=-\frac{\partial (K_1,\cdots,K_N)}{\partial (u_1,\cdots, u_N)}|_{\overline{u}}$
is negative definite with kernel space
$\{t\mathbf{1}\in \mathbb{R}^N|t\in \mathbb{R}\}$ generated by the normal vector of $\Sigma$ by Corollary \ref{Euclidean curvature Jacobian positivity}.
Note that the solution $u(t)$ of the extended modified Euclidean combinatorial Ricci flow (\ref{extended modified CRF equ}) stays in $\Sigma$.
This implies $\overline{u}$ is a local attractor of the extended modified Euclidean combinatorial Ricci flow (\ref{extended modified CRF equ}).
Then the exponential convergence of the solution $u(t)$ to $\overline{u}$
follows from the Lyapunov stability theorem (\cite{P}, Chapter 5).
\qed

\begin{remark}\label{extension not apply to CCF}
By Remark \ref{derivative tends infty Euclidean},
the extended combinatorial curvature $\widetilde{K}$ is not Lipschitz.
As a result, the combinatorial Laplace operator $\Delta=-\frac{\partial (K_1, \cdots, K_N)}{\partial (u_1, \cdots, u_N)}$
can not be extended by extending the combinatorial curvature $K$ to be $\widetilde{K}$.
Therefore, the combinatorial Calabi flow can not be extended in the way used for the combinatorial Ricci flow in this section.
In the special case that $\varepsilon\equiv0$, i.e. the case of Luo's vertex scaling,
there is another way introduced in \cite{GGLSW, GLSW} to extend the combinatorial Yamabe flow, where
one does surgery along the combinatorial Yamabe flow by edge flipping to preserve that the
triangulation is Delaunay along the combinatorial Yamabe flow.
The method of doing surgery by edge flipping under the Delaunay condition also applies to the combinatorial Calabi flow for Luo's vertex scalings \cite{ZX}.
It is proved that the solution of the combinatorial Yamabe flow with surgery \cite{GGLSW, GLSW}
and the combinatorial Calabi flow with surgery \cite{ZX} exists for all time and converges exponentially fast
for any initial piecewise linear and hyperbolic metric on the polyhedral surface.
In this approach, the assumption on the existence of discrete conformal factor with prescribed combinatorial curvature $\overline{K}$
is not needed. Please refer to \cite{GGLSW, GLSW, ZX} for more details on this.
In the special case that $\varepsilon\equiv 1$, i.e. the case of Bowers-Stephenson's inversive distance circle packings \cite{BS},
Bobenko-Lutz \cite{BL} recently developed similar surgery by edge flipping under the weighted Delaunay condition.
It is proved that the solution of the combinatorial Ricci flow with surgery \cite{BL}
and the combinatorial Calabi flow with surgery \cite{XZ2023} exists for all time and converges exponentially fast
for any initial decorated piecewise Euclidean metric on surfaces.
\end{remark}

\section{Hyperbolic combinatorial curvature flows}\label{section 3}

Paralleling to the Euclidean case, we study the following modified hyperbolic combinatorial Ricci flow
\begin{equation}\label{modified CRF equ hyper}
\begin{aligned}
\frac{du_i}{dt}=\overline{K}_i-K_i
\end{aligned}
\end{equation}
and modified hyperbolic combinatorial Calabi flow
\begin{equation}\label{modified CCF equ hyper}
\begin{aligned}
\frac{du_i}{dt}=\Delta (K-\overline{K})_i,
\end{aligned}
\end{equation}
where $\overline{K}: V\rightarrow (-\infty, 2\pi)$ is a function defined on the vertices
with $\sum_{i=1}^N\overline{K}_i>2\pi\chi(M)$ in the hyperbolic background geometry.
The modified hyperbolic combinatorial Ricci flow (\ref{modified CRF equ hyper})
and modified hyperbolic combinatorial Calabi flow (\ref{modified CCF equ hyper})
are generalizations of the combinatorial Ricci flow (\ref{CRF equ})
and the combinatorial Calabi flow (\ref{CCF equ}) respectively in the hyperbolic background geometry.
As many results in the hyperbolic case are paralleling to those in the Euclidean case,
we will just list the results and omit the proofs for the hyperbolic case if the proof has no difference from the Euclidean case.

To simplify the notations, set
\begin{equation}\label{simplification C_i S_i}
\begin{aligned}
S_i=e^{f_i}, C_i=\sqrt{1+\varepsilon_ie^{2f_i}}.
\end{aligned}
\end{equation}
Then
\begin{equation}\label{relation of C_i S_i}
\begin{aligned}
C_i^2-\varepsilon_iS_i^2=1
\end{aligned}
\end{equation}
and the edge length $l_{ij}$ is determined by
\begin{equation}\label{hyperbolic edge length in C_i S_i}
\begin{aligned}
\cosh l_{ij}=C_iC_j+\eta_{ij}S_iS_j.
\end{aligned}
\end{equation}

We recall some basic properties of Glickenstein's hyperbolic discrete conformal structures on surfaces.
Please refer to \cite{GT,Xu rigidity} for more details on these properties.

\begin{theorem}[\cite{Xu rigidity}]\label{hyperbolic admissible space}
Suppose $\sigma=\{ijk\}$ is a topological triangle with
two weights $\varepsilon: V_\sigma\rightarrow \{0, 1\}$ and $\eta: E_\sigma\rightarrow \mathbb{R}$
satisfying the structure conditions (\ref{structure condition 1}) and (\ref{structure condition 2}).
Then the admissible space $\Omega_{ijk}^H(\eta)\subseteq \mathbb{R}^3$ of
hyperbolic discrete conformal factors $f: V\rightarrow \mathbb{R}$ on $\{ijk\}$
is nonempty and simply connected with analytical boundary components.
Furthermore,
\begin{equation*}
\begin{aligned}
\Omega^H_{ijk}(\eta)=\mathbb{R}^3\setminus \sqcup_{\alpha\in \Lambda} V_\alpha,
\end{aligned}
\end{equation*}
where $\Lambda=\{q\in \{i,j,k\}|A_q=\eta_{st}^2-\varepsilon_s\varepsilon_t>0, \{q, s, t\}=\{i, j, k\}\}$,
$\sqcup_{\alpha\in \Lambda}V_\alpha$ is a disjoint union of closed regions in $\mathbb{R}^3$ bounded by an analytical function defined on $\mathbb{R}^2$
and $V_j, V_k$ defined similarly.
\end{theorem}

\begin{remark}[\cite{Xu rigidity}]\label{derivative tends infty hyperbolic}
If $(\eta_{ij}, \eta_{ik}, \eta_{jk})$ defined on the triangle $\{ijk\}$
satisfies the structure conditions (\ref{structure condition 1}) and (\ref{structure condition 2}) and
 $(f_i, f_j, f_k)\in \Omega_{ijk}^H(\eta)$ tends to a point $(\overline{f}_i, \overline{f}_j, \overline{f}_k)\in \partial V_i$ with $V_i\neq \emptyset$,
 then $\frac{\partial \theta_i}{\partial u_i}\rightarrow -\infty$, $\frac{\partial \theta_i}{\partial u_j}\rightarrow +\infty$, $\frac{\partial \theta_i}{\partial u_k}\rightarrow +\infty$.
\end{remark}

\begin{theorem}[\cite{GT, Xu rigidity}]\label{hyperbolic Jacobian negativity}
Suppose $\sigma=\{ijk\}$ is a topological triangle with
two weights $\varepsilon: V_\sigma\rightarrow \{0, 1\}$ and $\eta: E_\sigma\rightarrow \mathbb{R}$
satisfying the structure conditions (\ref{structure condition 1}) and (\ref{structure condition 2}).
Then the Jacobian matrix $\Lambda^H_{ijk}=\frac{\partial (\theta_i, \theta_j, \theta_k)}{\partial ( u_i, u_j, u_k)}$
is symmetric and strictly negative definite for all nondegenerate
hyperbolic discrete conformal factors on the triangle $\{ijk\}$.
\end{theorem}

\begin{corollary}[\cite{Xu rigidity}]\label{hyperbolic curvature Jacobian positivity}
Suppose $(M, \mathcal{T}, \varepsilon, \eta)$ is a weighted triangulated surface with
the weights $\varepsilon: V\rightarrow \{0, 1\}$ and $\eta: E\rightarrow \mathbb{R}$
satisfying the structure conditions (\ref{structure condition 1}) and (\ref{structure condition 2}).
Then the Jacobian matrix $\Lambda^H=\frac{\partial (K_1,\cdots,K_N)}{\partial (u_1,\cdots, u_N)}$
is symmetric and positive definite
for all nondegenerate hyperbolic discrete conformal factors on $(M, \mathcal{T}, \varepsilon, \eta)$.
\end{corollary}

Theorem \ref{hyperbolic admissible space} and Theorem \ref{hyperbolic Jacobian negativity} imply the following Ricci energy function for the triangle $\{ijk\}$
\begin{equation}\label{hyperbolic Ricci energy function for triangle}
\begin{aligned}
\mathcal{E}_{ijk}(u_i,u_j,u_k)=\int_{(\overline{u}_i, \overline{u}_j, \overline{u}_k)}^{(u_i, u_j, u_k)}\theta_idu_i+\theta_jdu_j+\theta_kdu_k
\end{aligned}
\end{equation}
is a well-defined smooth locally strictly concave function on $\Omega^H_{ijk}(\eta)$ with $\nabla_{u_i}\mathcal{E}_{ijk}=\theta_i$
under the structure conditions (\ref{structure condition 1}) and (\ref{structure condition 2}).
Set
\begin{equation}\label{hyperbolic Ricci energy function}
\begin{aligned}
\mathcal{E}(u_1, \cdots, u_N)=2\pi\sum_{i\in V}u_i-\sum_{\{ijk\}\in F}\mathcal{E}_{ijk}(u_i,u_j,u_k)
\end{aligned}
\end{equation}
to be the Ricci energy function defined on the admissible space $\Omega^H$ of nondegenerate hyperbolic discrete conformal factors
for $(M, \mathcal{T}, \varepsilon, \eta)$.
Then $\mathcal{E}$ is a locally strictly convex function defined on $\Omega^H$ with
$\nabla_{u_i} \mathcal{E}=K_i$
by Corollary \ref{hyperbolic curvature Jacobian positivity}.

Paralleling to Lemma \ref{CRF and CCF are gradient flows} in the Euclidean case,
we have the following result on the modified hyperbolic combinatorial Ricci flow (\ref{modified CRF equ  hyper}) and the modified hyperbolic combinatorial Calabi flow (\ref{modified CCF equ  hyper}).

\begin{lemma}\label{CRF and CCF are gradient flows hyperbolic}
Suppose $(M, \mathcal{T}, \varepsilon, \eta)$ is a weighted triangulated connected closed surface with
the weights $\varepsilon: V\rightarrow \{0, 1\}$ and $\eta: E\rightarrow \mathbb{R}$ satisfying
the structure conditions (\ref{structure condition 1}) and (\ref{structure condition 2}).
The modified hyperbolic combinatorial Ricci flow (\ref{modified CRF equ  hyper})
and modified hyperbolic combinatorial Calabi flow (\ref{modified CCF equ hyper}) on $(M, \mathcal{T}, \varepsilon, \eta)$
are negative gradient flows.
\end{lemma}

Paralleling to Theorem \ref{stability Euclidean} in the Euclidean case,
we have the following result on the longtime existence and convergence for the solutions of
the modified hyperbolic combinatorial Ricci flow (\ref{modified CRF equ hyper})
and the modified hyperbolic combinatorial Calabi flow (\ref{modified CCF equ hyper})
for initial values with small energy,
which is a slight generalization of Theorem \ref{convergence of CCF introduction} in the hyperbolic background geometry.

\begin{theorem}
Suppose $(M, \mathcal{T}, \varepsilon, \eta)$ is a weighted triangulated connected closed surface with
the weights $\varepsilon: V\rightarrow \{0, 1\}$ and $\eta: E\rightarrow \mathbb{R}$ satisfying
the structure conditions (\ref{structure condition 1}) and (\ref{structure condition 2}).
\begin{description}
  \item[(a)] If the solution of the modified hyperbolic combinatorial Ricci flow (\ref{modified CRF equ  hyper})
or modified hyperbolic combinatorial Calabi flow (\ref{modified CCF equ hyper}) converges to a nondegenerate hyperbolic discrete conformal
factor $\overline{u}$, then the combinatorial curvature for the hyperbolic polyhedral metric
determined by the discrete conformal factor $\overline{u}$ is $\overline{K}$.
  \item[(b)] If there exists a nondegenerate hyperbolic discrete conformal factor $\overline{u}$ with combinatorial curvature $\overline{K}$,
  then there exists a real number $\delta>0$ such that if the initial value $u(0)$ of
   the modified hyperbolic combinatorial Ricci flow (\ref{modified CRF equ  hyper}) (the modified hyperbolic combinatorial Calabi flow (\ref{modified CCF equ hyper}) respectively)
   satisfies $||u(0)-\overline{u}||<\delta$, the solution of the modified hyperbolic combinatorial Ricci flow (\ref{modified CRF equ hyper})
  (the modified hyperbolic combinatorial Calabi flow (\ref{modified CCF equ hyper}) respectively) exists for all time and converges exponentially fast to $\overline{u}$.
\end{description}
\end{theorem}

To handle the potential singularities that may develop along the modified hyperbolic combinatorial Ricci flow (\ref{modified CRF equ  hyper})
with general initial values,
we need the following extension introduced by the author \cite{Xu rigidity} for Glickenstein's hyperbolic discrete conformal structures.

\begin{lemma}\label{hyperbolic extension}
Suppose $\sigma=\{ijk\}$ is a topological triangle with
two weights $\varepsilon: V_\sigma\rightarrow \{0, 1\}$ and $\eta: E_\sigma\rightarrow \mathbb{R}$
satisfying the structure conditions (\ref{structure condition 1}) and (\ref{structure condition 2}).
Then the inner angle functions $\theta_i, \theta_j, \theta_k$ of the triangle $\{ijk\}$
defined for nondegenerate hyperbolic discrete conformal factors on $\{ijk\}$
can be extended by constants to be continuous functions $\widetilde{\theta}_i, \widetilde{\theta}_j, \widetilde{\theta}_k$
defined for $(f_i, f_j, f_k)\in \mathbb{R}^3$ by setting
\begin{equation}\label{extension of theta_i hyperbolic}
\begin{aligned}
\widetilde{\theta}_i(f_i,f_j,f_k)=\left\{
                       \begin{array}{ll}
                         \theta_i, & \hbox{if $(f_i,f_j,f_k)\in \Omega^H_{ijk}(\eta)$;} \\
                         \pi, & \hbox{if $(f_i,f_j,f_k)\in V_i$;} \\
                         0, & \hbox{otherwise.}
                       \end{array}
                     \right.
\end{aligned}
\end{equation}
\end{lemma}

By the extension in Lemma \ref{hyperbolic extension}, the combinatorial curvature $K$ defined for nondegenerate hyperbolic discrete conformal factors in $\Omega^H$ can be extended to be the generalized combinatorial curvature $\widetilde{K}_i=2\pi-\sum\widetilde{\theta}_i$  defined for any $f\in \mathbb{R}^N$.
We call $f\in \mathbb{R}^N$ as a generalized hyperbolic discrete conformal factor.

Taking $\widetilde{\theta}_i, \widetilde{\theta}_j, \widetilde{\theta}_k$ as functions of $(u_i, u_j, u_k)$.
Then the extensions $\widetilde{\theta}_i, \widetilde{\theta}_j, \widetilde{\theta}_k$ of $\theta_i, \theta_j, \theta_k$
 are continuous functions of
$(u_i, u_j, u_k)\in V_i\times V_j\times V_k$, where $V_q=\mathbb{R}$ if $\varepsilon_q=0$ and
$V_q=\mathbb{R}_{<0}=(-\infty, 0)$ if $\varepsilon_q=1$ for $q\in \{i,j,k\}$.
Combining this with Theorem \ref{Luo's convex extention}, the locally concave function
$\mathcal{E}_{ijk}$ defined by (\ref{hyperbolic Ricci energy function for triangle})
for nondegenerate hyperbolic discrete conformal factors for a triangle $\{ijk\}$
can be extended to be a $C^1$ smooth concave function
\begin{equation}\label{extension of hyperbolic Ricci energy function for triangle}
\begin{aligned}
\widetilde{\mathcal{E}}_{ijk}(u_i,u_j,u_k)=\int_{(\overline{u}_i, \overline{u}_j, \overline{u}_k)}^{(u_i, u_j, u_k)}\widetilde{\theta}_idu_i+\widetilde{\theta}_jdu_j+\widetilde{\theta}_kdu_k
\end{aligned}
\end{equation}
defined on $V_i\times V_j\times V_k$ with $\nabla_{u_i}\widetilde{\mathcal{E}}_{ijk}=\widetilde{\theta}_i$.
As a result, the locally convex function $\mathcal{E}$ defined by (\ref{hyperbolic Ricci energy function}) for nondegenerate hyperbolic
discrete conformal factors on a weighted triangulated surface $(M, \mathcal{T}, \varepsilon, \eta)$ can be extended to
be a $C^1$ smooth convex function
\begin{equation}\label{extended hyperbolic energy function}
\begin{aligned}
\widetilde{\mathcal{E}}(u_1, \cdots, u_N)=2\pi\sum_{i\in V}u_i-\sum_{\{ijk\}\in F}\widetilde{\mathcal{E}}_{ijk}(u_i,u_j,u_k)
\end{aligned}
\end{equation}
defined on $\mathbb{R}^{N_1}\times \mathbb{R}^{N_2}_{<0}$ with $\nabla_{u_i} \widetilde{\mathcal{E}}=\widetilde{K}_i=2\pi-\sum\widetilde{\theta}_i$,
where $N_1$ is the number of vertices $v_i$ in $V$ with $\varepsilon_i=0$ and $N_2=N-N_1$.
We also call a vector $u\in \mathbb{R}^{N_1}\times \mathbb{R}^{N_2}_{<0}$ as
a generalized hyperbolic discrete conformal factor for the
weighted triangulated surface $(M, \mathcal{T}, \varepsilon, \eta)$.

Using the function $\widetilde{\mathcal{E}}$ defined by (\ref{extended hyperbolic energy function}), the author \cite{Xu rigidity} proved the following rigidity for Glickenstein's hyperbolic discrete conformal structures.

\begin{theorem}\label{Thm rigidity of HDCS context}
Suppose $(M, \mathcal{T}, \varepsilon, \eta)$ is a weighted triangulated surface with
the weights $\varepsilon: V\rightarrow \{0, 1\}$ and $\eta: E\rightarrow \mathbb{R}$
satisfying the structure conditions (\ref{structure condition 1}) and (\ref{structure condition 2}).
If there exists a nondegenrate hyperbolic discrete conformal factor $f_A\in \Omega^H$ and
a generalized hyperbolic discrete conformal factor $f_B\in \mathbb{R}^N$ such that
$K(f_A)=\widetilde{K}(f_B)$. Then $f_A=f_B$.
\end{theorem}

Using the extension of inner angles in Lemma \ref{hyperbolic extension}, we can extend the modified hyperbolic combinatorial Ricci flow (\ref{modified CRF equ hyper}) to the following form.

\begin{definition}
Suppose $(M, \mathcal{T}, \varepsilon, \eta)$ is a weighted triangulated connected closed surface with
the weights $\varepsilon: V\rightarrow \{0, 1\}$ and $\eta: E\rightarrow \mathbb{R}$ satisfying
the structure conditions (\ref{structure condition 1}) and (\ref{structure condition 2}).
The extended modified hyperbolic combinatorial Ricci flow is defined to be
\begin{equation}\label{extended modified CRF equ hyper}
\begin{aligned}
\frac{du_i}{dt}=\overline{K}_i-\widetilde{K}_i,
\end{aligned}
\end{equation}
where $\widetilde{K}_i=2\pi-\sum_{\{ijk\}\in F}\widetilde{\theta}_i$ is an extension of the combinatorial curvature $K_i$ with
$\widetilde{\theta}_i$ given by  (\ref{extension of theta_i hyperbolic}) in the hyperbolic background geometry.
\end{definition}

Paralleling to Theorem \ref{Euclidean uniqueness} in the Euclidean case, we have the following uniqueness for the solution of the
extended modified hyperbolic combinatorial Ricci flow (\ref{extended modified CRF equ hyper}).
\begin{theorem}\label{hyperbolic uniqueness}
Suppose $(M, \mathcal{T}, \varepsilon, \eta)$ is a weighted triangulated connected closed surface with
the weights $\varepsilon: V\rightarrow \{0, 1\}$ and $\eta: E\rightarrow \mathbb{R}$ satisfying
the structure conditions (\ref{structure condition 1}) and (\ref{structure condition 2}).
Then the solution of the extended modified hyperbolic combinatorial Ricci flow (\ref{extended modified CRF equ hyper})
is unique for any initial generalized
hyperbolic discrete conformal factor $f(0)\in\mathbb{R}^N$ on $(M, \mathcal{T}, \varepsilon, \eta)$.
\end{theorem}

To prove the hyperbolic version of Theorem \ref{longtime existence and convgence of CRF Euclidean},
we need some more preparations.

\begin{lemma}\label{estimate of cosh l using eta}
Suppose $i, j$ are two adjacent vertices in $V$ and the weight $\eta_{ij}$ satisfies the structure condition (\ref{structure condition 1}).
If the edge length $l_{ij}$ is defined by (\ref{defn of hyperbolic length}), $\varepsilon_i=1$ and $\varepsilon_j\in \{0,1\}$,
then there exist two positive constants $\lambda=\lambda(\varepsilon_j, \eta_{ij})$ and $\mu=\mu(\eta_{ij})$ such that
\begin{equation*}
\begin{aligned}
\lambda(C_iC_j+S_iS_j)\leq \cosh l_{ij}\leq \mu(C_iC_j+S_iS_j),
\end{aligned}
\end{equation*}
where $C_i, C_j, S_i, S_j$ are defined by (\ref{simplification C_i S_i}).
\end{lemma}
\proof
By (\ref{defn of hyperbolic length}), the edge length $l_{ij}$ satisfies
$\cosh l_{ij}=C_iC_j+\eta_{ij}S_iS_j\leq (1+|\eta_{ij}|)(C_iC_j+S_iS_j)$.
Therefore, we can take $\mu=1+|\eta_{ij}|$.

If $\varepsilon_j=1$, then $\eta_{ij}>-1$ by the structure condition (\ref{structure condition 1}).
In this case, $C_i=\sqrt{1+e^{2f_i}}>e^{f_i}=S_i$ and similarly $C_j>S_j$.
If $\eta_{ij}>0$, then $\cosh l_{ij}\geq\min \{1, \eta_{ij}\}(C_iC_j+S_iS_j)$, where $\min \{1, \eta_{ij}\}>0$.
If $-1<\eta_{ij}\leq 0$, by $C_i>S_i, C_j>S_j$, we have
$\cosh l_{ij}\geq (1+\eta_{ij})C_iC_j>\frac{1}{2}(1+\eta_{ij})(C_iC_j+S_iS_j)$. We can take $\lambda=\frac{1}{2}(1+\eta_{ij})>0$ in this case.

If $\varepsilon_j=0$, then $\eta_{ij}>0$ by the structure condition (\ref{structure condition 1}).
In this case, we have $\cosh l_{ij}\geq \min \{1, \eta_{ij}\}(C_iC_j+S_iS_j)$ with $\min \{1, \eta_{ij}\}>0$.
Therefore, we can take $\lambda=\min \{1, \eta_{ij}\}>0$ in this case.
\qed

\begin{lemma}\label{small theta i for large fi}
Suppose $\sigma=\{ijk\}$ is a topological triangle with
two weights $\varepsilon: V_\sigma\rightarrow \{0, 1\}$ and $\eta: E_\sigma\rightarrow \mathbb{R}$
satisfying $\varepsilon_i=1$ and the structure conditions (\ref{structure condition 1}) and (\ref{structure condition 2}).
Then for any $\epsilon>0$, there exists
a positive number $L=L(\varepsilon,\eta,\epsilon)$ such that if $f_i>L$, the inner angle $\theta_i$ at the vertex $i$
of the nondegenerate hyperbolic triangle $\{ijk\}$ with edge lengths defined by (\ref{defn of hyperbolic length})
is smaller than $\epsilon$.
\end{lemma}
\proof
By the hyperbolic cosine law, we have
\begin{equation}\label{estimate of cosine theta in nu and omega}
\begin{aligned}
\cos \theta_i
=&\frac{\cosh l_{ij}\cosh l_{ik}-\cosh l_{jk}}{\sinh l_{ij}\sinh l_{ik}}\\
=&\frac{\cosh(l_{ij}+l_{ik})+\cosh (l_{ij}-l_{ik})-2\cosh l_{jk}}{\cosh(l_{ij}+l_{ik})-\cosh (l_{ij}-l_{ik})}\\
=&\frac{1+\nu-2\omega}{1-\nu},
\end{aligned}
\end{equation}
where $\nu=\frac{\cosh (l_{ij}-l_{ik})}{\cosh(l_{ij}+l_{ik})}$ and $\omega=\frac{\cosh l_{jk}}{\cosh(l_{ij}+l_{ik})}$.
By Lemma \ref{estimate of cosh l using eta}, we have
\begin{equation}\label{estimate of nu}
\begin{aligned}
0<\nu<\frac{\cosh (l_{ik})}{\cosh(l_{ij}+l_{ik})}<\frac{1}{\cosh l_{ij}}\leq\frac{1}{\lambda(C_iC_j+S_iS_j)}<\frac{1}{\lambda C_i}<\frac{1}{\lambda S_i}
\end{aligned}
\end{equation}
and
\begin{equation}\label{estimate of omega}
\begin{aligned}
0<\omega<\frac{\mu(C_jC_k+S_jS_k)}{\lambda^2(C_iC_j+S_iS_j)(C_iC_k+S_iS_k)}
<\frac{\mu(C_jC_k+S_jS_k)}{\lambda^2(C_i^2C_jC_k+S_i^2S_jS_k)}<\frac{\mu}{\lambda^2S_i^2},
\end{aligned}
\end{equation}
where $C_i>S_i$ is used for $\varepsilon_i=1$.
Note that $S_i=e^{f_i}$, (\ref{estimate of nu}) and (\ref{estimate of omega}) imply $\nu,\omega\rightarrow 0$ uniformly as $f_i\rightarrow +\infty$.
By (\ref{estimate of cosine theta in nu and omega}), $\theta_i$ tends to $0$ uniformly as $f_i\rightarrow +\infty$.
Therefore, for any $\epsilon>0$, there exists $L>0$ such that if $f_i>L$, then $\theta_i<\epsilon$.
\qed


As a corollary of Lemma \ref{small theta i for large fi}, we have
the following result.

\begin{corollary}\label{small extended theta i for large fi}
Suppose $\sigma=\{ijk\}$ is a topological triangle with
two weights $\varepsilon: V_\sigma\rightarrow \{0, 1\}$ and $\eta: E_\sigma\rightarrow \mathbb{R}$
satisfying $\varepsilon_i=1$ and the structure conditions (\ref{structure condition 1}) and (\ref{structure condition 2}).
Then for any $\epsilon>0$, there exists
a positive number $L=L(\varepsilon,\eta,\epsilon)$ such that if $f_i>L$, the extended inner angle $\widetilde{\theta}_i$
defined by (\ref{extension of theta_i hyperbolic})
at the vertex $i$
in the generalized hyperbolic triangle $\{ijk\}$ with edge lengths defined by (\ref{defn of hyperbolic length})
is smaller than $\epsilon$.
\end{corollary}
\proof
If the generalized hyperbolic triangle does not degenerate, then by Lemma \ref{small theta i for large fi},
for $\epsilon>0$, there exists a constant $L=L(\varepsilon,\eta,\epsilon)>0$ such that if $f_i>L$, then $\widetilde{\theta}_i=\theta_i<\epsilon$.
If the generalized hyperbolic triangle $\{ijk\}$ degenerates, we claim that if $f_i>L$, then $\widetilde{\theta}_i=0<\epsilon$.
Otherwise, $\widetilde{\theta}_i=\pi$ by Lemma \ref{hyperbolic extension}.
Combining the continuity of $\widetilde{\theta}_i$ and the fact that $\widetilde{\theta}_i$ is a constant extension of $\theta_i$, we have
$\widetilde{\theta}_i\leq \epsilon<\pi$. It is a contradiction.
\qed

Now we can prove the hyperbolic version of Theorem \ref{longtime existence and convgence of CRF Euclidean},
which generalizes Theorem \ref{convergence of extended CRF introduction} (a) (c) in the hyperbolic background geometry.
\begin{theorem}\label{longtime existence and convgence of CRF hyperbolic}
Suppose $(M, \mathcal{T}, \varepsilon, \eta)$ is a weighted triangulated connected closed surface with
the weights $\varepsilon: V\rightarrow \{0, 1\}$ and $\eta: E\rightarrow \mathbb{R}$ satisfying
the structure conditions (\ref{structure condition 1}) and (\ref{structure condition 2}).
Then the solution of the extended modified hyperbolic combinatorial Ricci flow (\ref{extended modified CRF equ hyper}) exists for all time
for any initial generalized hyperbolic discrete conformal factor $u$ on $(M, \mathcal{T}, \varepsilon, \eta)$.
Furthermore, if there exists a nondegenerate hyperbolic discrete conformal factor $\overline{u}$ with combinatorial curvature $\overline{K}$,
then the solution of the extended modified hyperbolic combinatorial Ricci flow (\ref{extended modified CRF equ hyper})
converges exponentially fast to $\overline{u}$ for any initial generalized hyperbolic discrete conformal factor $u(0)$.
\end{theorem}
\proof
By (\ref{u hyperbolic}), $u_i\in \mathbb{R}$ for any vertex $i$ with $\varepsilon_i=0$ and $u_i\in \mathbb{R}_{<0}=(-\infty, 0)$ for any vertex $i$ with $\varepsilon_i=1$. Therefore, $u=(u_1,\cdots, u_N)\in \mathbb{R}^{N_1}\times \mathbb{R}_{<0}^{N_2}$,
where $N_1$ is the number of vertices $i\in V$ with $\varepsilon_i=0$ and $N_2=N-N_1$.
If $u(t)$ is a solution of the extended modified hyperbolic combinatorial Ricci flow (\ref{extended modified CRF equ hyper}),
then $|u_i(t)|\leq |u_i(0)|+[|\overline{K}_i|+(d_i+2)\pi]t<+\infty$ for $t\in [0, +\infty)$,
where $d_i$ is the number of vertices adjacent to the vertex $i\in V$.
This implies $u_{i}(t)$ is bounded for the vertex $i\in V$ with $\varepsilon_i=0$ and bounded from below
for $i\in V$ with $\varepsilon_i=1$ in finite time.

We claim that $u_i(t)$ is uniformly bounded from above in $(-\infty, 0)$ for $i\in V$ with $\varepsilon_i=1$.
Then the longtime existence for the solution of the extended modified hyperbolic combinatorial Ricci flow (\ref{extended modified CRF equ}) follows.
We use Ge-Xu's trick in \cite{GX4} to prove the claim by contradiction.
Suppose there exists a vertex $i\in V$ such that $\overline{\lim}_{t\uparrow T} u_i(t)=0$
for $T\in(0, +\infty]$,
which corresponds to $\overline{\lim}_{t\uparrow T} f_i(t)=+\infty$ by (\ref{u hyperbolic}).
By Corollary \ref{small extended theta i for large fi}, for $\epsilon=\frac{1}{d_i}(2\pi-\overline{K})>0$,
where $d_i$ is the degree of the vertex $v_i$,
there exists a constant $c<0$ such that if $u_i>c$, then $\widetilde{\theta}_i<\epsilon$ and hence $\widetilde{K}_i>\overline{K}_i$.
Choose $t_0\in (0,T)$ such that $u_i(t_0)>c$, the existence of which is ensured by $\overline{\lim}_{t\uparrow T} u_i(t)=0$.
Set $a=\inf \{t<t_0|u_i(s)>c, \forall s\in [t, t_0]\}$, then $u_i(a)=c$.
Note that for $t\in (a, t_0]$, $u_i'(t)=\overline{K}-\widetilde{K}_i<0$ along the flow (\ref{extended modified CRF equ}),
we have $u_i(t_0)<u(a)=c$, which contradicts the assumption that $u_i(t_0)>c$.
The arguments here further imply that $u_i(t)$ is uniformly bounded from above in $(-\infty, 0)$ for all $i\in V$ with $\varepsilon_i=1$.

If there exists a nondegenerate hyperbolic discrete conformal factor $\overline{u}$ with combinatorial curvature $\overline{K}$,
then $\overline{u}$ is a critical point of the $C^1$ smooth convex function
$\widetilde{\mathcal{H}}(u)=\widetilde{\mathcal{E}}(u)-\int_{\overline{u}}^u\sum_{i=1}^N\overline{K}_idu_i$, where $\widetilde{\mathcal{E}}(u)$
is the extended Ricci energy function defined by (\ref{extended hyperbolic energy function}).
Note that $0=\widetilde{\mathcal{H}}(\overline{u})\leq \widetilde{\mathcal{H}}(u)$ and $\nabla \widetilde{\mathcal{H}}(\overline{u})=0$,
by Lemma \ref{convex func tends infty at infty},
we have $\lim_{u\rightarrow \infty}\widetilde{\mathcal{H}}(u)=+\infty$.
Further note that
\begin{equation*}
\frac{d}{dt}\widetilde{\mathcal{H}}(u(t))=\nabla_u \widetilde{\mathcal{H}}\cdot \frac{du}{dt}=-\sum_{i=1}^N(\widetilde{K}_i-\overline{K}_i)^2\leq 0
\end{equation*}
along the extended modified hyperbolic combinatorial Ricci flow (\ref{extended modified CRF equ hyper}),
we have $0\leq \widetilde{\mathcal{H}}(u(t))\leq \widetilde{\mathcal{H}}(u(0))$.
This implies that the solution $u(t)$ of the extended modified hyperbolic combinatorial Ricci flow (\ref{extended modified CRF equ hyper})
lies in a compact subset of $\mathbb{R}^N$.
Combining with the fact that  $u_i(t)$ is uniformly bounded from above in $(-\infty, 0)$ for any vertex $i\in V$ with $\varepsilon_i=1$,
the solution $u(t)$ of the extended modified hyperbolic combinatorial Ricci flow (\ref{extended modified CRF equ hyper})
lies in a compact subset of $\mathbb{R}^{N_1}\times \mathbb{R}_{<0}^{N_2}$.
The proof in the following is the same as that for Theorem \ref{longtime existence and convgence of CRF Euclidean},
so we omit the details here.
\qed

\section{Open problems}\label{section 4}

\subsection{Convergence of combinatorial curvature flows with surgery}

In Theorem \ref{convergence of extended CRF introduction}, we extend the combinatorial Ricci flow
through the singularities of the flow to ensure the convergence of the flow under the assumption that
there exists a discrete conformal factor with constant combinatorial curvature.
On the one hand, this method can not be applied to the combinatorial Calabi flow by
Remark \ref{derivative tends infty Euclidean} and Remark \ref{derivative tends infty hyperbolic}.
On the other hand, we do not hope the combinatorial curvature flows develop singularities in practical applications.
One way to avoid the singularities is to do surgery along the combinatorial curvature flows before the singularities develops.
Motivated by the surgery by edge flipping under the Delaunay condition introduced in \cite{GLSW,GGLSW} for Luo's vertex scaling,
and the the surgery by edge flipping under the weighted Delaunay condition introduced in \cite{BL, BL2} for decorated piecewise Euclidean metrics,
we can also do surgery by edge flipping under the weighted Delaunay condition along combinatorial curvature flows for Glickenstein's discrete conformal structures on polyhedral surfaces. Motivated by Glickenstein's work \cite{G2a,G3,G4}, we introduce the following definition of weighted Delaunay triangulation for Glickenstein's discrete conformal structures on polyhedral surfaces.

Suppose $(M, V)$ is a marked surface and $V$ is a nonempty finite subset of $M$.
$\varepsilon: V\rightarrow \{0,1\}$ is a weight defined on $V$. $(M, V, \varepsilon)$ is called a weighted marked surface.
\begin{definition}[\cite{G2a,G3,G4}]
Suppose $(M, V, \varepsilon)$ is a weighted marked surface with a Euclidean polyhedral metric $d$ and
$\mathcal{T}$ is a geometric triangulation of $(M, V, \varepsilon, d)$ with every triangle $\{ijk\}$  in the triangulation have a
well-defined geometric center $C_{ijk}$.
Suppose $\{ij\}$ is an edge shared by two adjacent Euclidean triangles $\{ijk\}$ and $\{ijl\}$.
The edge $\{ij\}$ is called weighted Delaunay if $h_{ij,k}+h_{ij,l}\geq 0$, where
$h_{ij,k}, h_{ij,l}$ are the signed distance of $C_{ijk}, C_{ijl}$ to the edge $\{ij\}$ respectively.
The triangulation $\mathcal{T}$ is called weighted Delaunay in $d$ if every edge in the triangulation is weighted Delaunay.
\end{definition}

Along the Euclidean combinatorial curvature flows (the Eucldiean combinatorial Ricci flow or the combinatorial Calabi flow)
for Glickenstein's Euclidean discrete conformal structures on a weighted triangulated surface $(M, V, \varepsilon, \mathcal{T})$,
if the triangulation $\mathcal{T}$ is weighted Delaunay in the induced Euclidean polyhedral metric $d(u(t))$  for
$t\in [0, T]$ and not weighted Delaunay in $d(u(t))$  for $t\in (T, T+\epsilon)$, $\epsilon>0$, there exists
an edge $\{ij\}\in E$ such that $h_{ij,k}+h_{ij,l}\geq 0$
for $t\in [0, T]$ and $h_{ij,k}+h_{ij,l}<0$ for $t\in (T, T+\epsilon)$.
We replace the triangulation $\mathcal{T}$ by a new triangulation $\mathcal{T}'$ at the time $t=T$
by replacing two triangles $\{ijk\}$ and $\{ijl\}$
adjacent to $\{ij\}$ by two new triangles $\{ikl\}$ and $\{jkl\}$.
This is called a \textbf{surgery by flipping} under the weighted Delaunay condition on the triangulation $\mathcal{T}$,
which is an isometry of $(M, V, \varepsilon)$ in the Euclidean polyhedral metric $d(u(T))$.
After the surgery at time $t=T$, we run the Euclidean combinatorial curvature flow
on $(M, V,\varepsilon, \mathcal{T}')$ with initial metric coming from
the Euclidean combinatorial curvature flow on $(M, V, \varepsilon, \mathcal{T})$ at time $t=T$.
The surgery by flipping for hyperbolic combinatorial curvature flows can defined similarly.

We have the following conjecture on the longtime existence and convergence for the combinatorial
Ricci flow and the combinatorial Calabi flow with surgery, which is supported by the results in \cite{GLSW,GGLSW,ZX,XZ2023,XZ2023b,BL,BL2}
for some special cases.

\begin{conjecture}
Suppose $(M, V, \varepsilon)$ is a closed connected weighted marked surface with $\varepsilon: V\rightarrow \{0,1\}$.
For any initial Euclidean or hyperbolic polyhedral metric induced by Glickenstein's discrete conformal structures on $(M, V, \varepsilon)$,
the solution of the combinatorial Ricci flow and the combinatorial Calabi flow with surgery
exists for all time and converges exponentially fast.
\end{conjecture}


%
%
%
%

\end{document}